\documentclass{elsart}

\def\aa{\mathop{\bf a}\nolimits}
\def\x{\mathop{\bf x}\nolimits}
\def\y{\mathop{\bf y}\nolimits}
\def\z{\mathop{\bf z}\nolimits}
\def\a{\mathop{\alpha}\nolimits}
\def\b{\mathop{\beta}\nolimits}

\def\ax{a(\x)}
\def\bx{b(\x)}
\def\ux{u(\x)}
\def\fx{f(\x)}
\def\fhat{\widehat{f}}

\def\_#1{\mathop{\!^{}_{#1}}}

\def\ok{1,\ldots,k}
\def\fP{f\hspace{-.4ex}\_P}
\def\fPx{\fP(\x)}
\def\R{\mathop{{\rm I\!R}}\nolimits}
\def\Rk{\R^k}
\def\tR{\widetilde{\R}}

\begin{document}

\begin{frontmatter}
\title{Extending Utility Representations\\ of Partial Orders}
\date{April 8, 2001}

\author{Pavel Chebotarev\thanksref{coord}}
\address{Trapeznikov Institute of Control and Management Sciences,
65 Profsoyuznaya, Moscow 117997, Russia
}
\thanks[coord]{%
Fax: +7-095-420-20-16.\\
{\em E-mail addresses:} chv@lpi.ru, pchv@rambler.ru}

\begin{abstract}
The problem is considered as to whether a monotone function defined on a
subset $P$ of Euclidean space $\Rk$ can be strictly monotonically
extended to $\Rk$. It is proved that this is the case if and only if
the function is {\em separably increasing}. Explicit formulas are given
for a class of extensions which involves an arbitrary function. Similar
results are obtained for utility functions that represent strict
partial orders on abstract sets~$X$. The special case where $P$ is a
Pareto subset of~$\Rk$ (or of~$X$) is considered.\\
\end{abstract}

\end{frontmatter}

\section{Introduction}

Suppose that a decision maker defines his or her utility function on
some subset $P$ of the Euclidean space of alternatives~$\Rk$.
Utility functions are usually assumed to be strictly increasing
in the coordinates which correspond to partial criteria. Therefore, it
is useful to specify the conditions under which a strictly monotone
function defined on $P$ can be strictly monotonically extended
to~$\Rk$.

In this paper, we demonstrate that such an extension is possible if and
only if the function defined on $P$ is {\em separably increasing}.
Explicit formulas for the extension are provided.

An interesting special case is where the structure of subset $P$
does not permit any violation of strict monotonicity on~$P$. This
is the case where $P$ is a {\em Pareto set}. A corollary given in
Section~\ref{SPar} addresses this situation.
The results are translated to the general case of utility functions
that represent strict partial orders on arbitrary sets. We do not
impose continuity requirements in these settings.

\section{Notation, definitions, and main results}

For any $\x,\y\in \Rk,\:$
$\x\ge \y$ means [$x\_i\ge y\_i$ for all $i\in\{\ok\}$];
$\x\le \y$ means [$x\_i\le y\_i$ for all $i\in\{\ok\}$];
$\x >  \y$ means [$\x\ge \y$ and not $\x=\y$];
$\x <  \y$ means [$\x\le \y$ and not $\x=\y$].
These relations $\ge$, $\le$, $>$, and $<$ on $\Rk$ will be called
{\em Paretian}.

Consider an arbitrary subset $P$ of $\Rk$ and strictly increasing
(with respect to the above $>$ relation) real-valued functions $\fPx$
defined on $P$. The problem is to monotonically extend such a function
$\fPx$ to $\Rk,$ provided that it is possible, and to indicate
conditions under which this is possible.

An arbitrary function $\fP(\x)$ defined on any $P\subseteq\Rk$ is said
to be {\em strictly increasing on $P$ with respect to the $>$
relation\/} or simply {\em strictly increasing on $P$\/} if for every
$\x,\y\in P,\,$ $\x>\y$ implies $\fPx>\fP(\y)$.

\begin{defn}
\label{def2}
A real-valued function $\fPx$ defined on $P\subset \Rk$ is {\em
monotonically extendible to\/ $\Rk$} if there exists a function
$\fx:\;\Rk\to\R$ such that
\\$(\ast)$ the restriction of $\fx$ to $P$ coincides with
$\fPx,$ and
\\$(\ast\ast)$ $\fx$ is strictly increasing with respect to~$>$.\\
In this case$,$ $\fx$ is a {\em monotone extension\/} of $\fPx$ to~$\Rk$.
\end{defn}

The functions $\fPx$ and $\fx$ can be naturally considered as utility
(or objective) functions. This means that $\fPx$ and $\fx$ can be
interpreted as real-valued functions that represent the preferences of
a decision maker on the corresponding sets of alternatives (the
alternatives are identified with $k$-dimensional vectors of partial
criteria values or vectors of goods).

The Paretian $>$ relation is a strict partial order on~$\Rk$, i.e., it
is transitive and irreflexive. That is why we deal with a specific
problem on the monotonic extension of functions which are monotone with
respect to a partial order on their domain. We discuss some connections
of this problem with the classical results of utility theory in
Section~\ref{SUti} and turn to a more general formulation in
Section~\ref{abstSet}.

For technical convenience, let us add two extreme points to $\R$:
$$
\tR=\R\cup\,\{-\infty,+\infty\},
$$
and extend the ordinary $>$ relation to $\tR$: for every $\x\in \R,$ set
$+\infty>\x>-\infty$ and $+\infty>-\infty$. This $>$ relation will
determine the values of $\min$ and $\max$ functions on finite subsets
of~$\tR$.

The functions $\,\sup\,$ and $\,\inf\,$ will be considered as maps of
$2^{\R}$ to $\tR$ defined for the empty set as follows:
$\sup\emptyset=-\infty$ and
$\inf\emptyset=+\infty$.

Given $\fPx$, define two auxiliary functions for every $\x\in\Rk$:
\begin{eqnarray}
\label{ax}
\ax&=\sup\big\{\fP(\y)\mid\y\le\x,&\y\in P\big\},\\
\label{bx}
\bx&=\inf\big\{\fP(\z)\mid\z\ge\x,&\z\in P\big\}.
\end{eqnarray}
By this definition, $\ax$ and $\bx$ can be infinite.

It follows from the transitivity of the Paretian $>$ relation that for
every function $\fPx,$ functions $\ax$ and $\bx$ are nonstrictly
increasing with respect to~$>$:
\begin{equation}
\label{abmono}
\mbox{For all } \x,\x'\in \Rk,\,
\x'>\x \mbox{ implies } \big[a(\x')\ge\ax\mbox{ and }b(\x')\ge\bx\big].
\end{equation}
Moreover,
\begin{equation}
\mbox{For any }\,\x\in P,\quad\ax\ge\fPx\ge\bx.
\label{a>f>b}
\end{equation}

It can be easily shown that $\fPx$ is nonstrictly increasing if and
only if \begin{equation} \bx\ge\ax\quad\mbox{ for all }\x\in\Rk.
\label{b>a}
\end{equation}

The following definition involves a strengthening of~(\ref{b>a}).

\begin{defn}
\label{def1}
A~function $\fPx$ defined on $P\subset \Rk$ is {\em separably
increasing\/} if for any $\x,\x'\in\Rk,\:$ $\x'>\x$ implies
$b(\x')>\ax$.
\end{defn}

\begin{defn}
\label{def23}
Given $P\subseteq \Rk,$ let us say that
$P'\subseteq P$ is an {\em upper set\/} if for some $\aa\in\Rk,\,$
$P' =\big\{\x\mid\x\ge\aa,\,\x\in P\big\};$
$P''\subseteq P$ is a {\em lower set\/} if for some $\aa\in\Rk,\,$
$P''=\big\{\x\mid\x\le\aa,\,\x\in P\big\}$.
\end{defn}

\begin{prop}
\label{P0}
If $\fPx$ defined on $P\subset \Rk$ is separably increasing$,$ then\\
$(a)$ $\fPx$ is strictly increasing$;$\\
$(b)$ $\fPx$ is upper-bounded on lower sets and lower-bounded on upper
sets$;$ in other terms$,$ there are no $\x\in\Rk$ such that
$a(\x)=+\infty$ or $\/b(\x)=-\infty;$\\
$(c)$ For every\/ $\x\in\Rk,\:\bx\ge\ax;$\\
$(d)$ For every\/ $\x\in P,\:\bx=\ax=\fPx.$
\end{prop}

All proofs are given in Section~\ref{proofs}. Proposition~\ref{P0} and
other statements are proved there in the more general case of
utility functions that represent strict partial orders on arbitrary
sets (cf.\ Section~\ref{abstSet}).

Observe that there are functions $\fPx$ that are strictly increasing,
upper-bounded on lower sets and lower-bounded on upper sets, but are
not separably increasing. An example is
\begin{equation}
\fPx=\cases{x\_1,   &$x\_1\le0,$\cr
            x\_1-1, &$x\_1>1,$\cr
}
\label{CoEx}
\end{equation}
where $P=\,]-\infty,0]\,\cup\,]1,+\infty]\subset\R^1.$ This function
satisfies (a) and (b) (and, as well as all nonstrictly increasing
functions, it also satisfies (c) and (d)) of Proposition~\ref{P0}, but
it is not separably increasing.  Indeed, $b(1)=0=a(0).$

Suppose that $\fPx$ defined on $P\subset\Rk$ is separably increasing.
Below we prove that this is a necessary and sufficient condition for
the existence of monotone extensions of $\fPx$ to~$\Rk$ and demonstrate
how such extensions can be constructed.

Let $\ux\!:\Rk\to\R$ be any strictly increasing (with respect to the
Paretian $>$) and bounded function defined on the whole space~$\Rk$.
Suppose that $\a,\b\in\R$ are such that
\begin{equation}
\a<\ux<\b\quad\mbox{for all }\x\in\Rk. \label{ge}
\end{equation}

As an example of such a function, adduce
\begin{equation}
\label{arctan}
u\_{{\rm
example}}(\x)=
\frac{\b-\a}{\pi}\left(\arctan\sum_{i=1}^kx\_i+\frac{\pi}{2}\right)+\a.
\end{equation}

Consider also the special  case of strictly increasing functions
$u\_1(\x)$ such that
\begin{equation}
0<u\_1(\x)<1.
\label{0<g1<1}
\end{equation}
These functions can be obtained from the strictly increasing functions
$u(\x)$ that satisfy (\ref{ge}) as follows:
\begin{equation}
\label{u1}
u\_1(\x)=(\b-\a)^{-1}(\ux-\a).
\end{equation}

For an arbitrary strictly increasing function $u\_1(\x)$ that
satisfies (\ref{0<g1<1}), every $\a$ and $\b>\a,$ and every
$\x\in\Rk$, let us define
\begin{eqnarray}
\fx&=&
\max\Big\{\ax,\,\min\big\{\bx,\,\b\big\}-\b+\a\!\Big\}\big(1-u\_1(\x)\big)
\nonumber\\
   &+&\min\Big\{\bx,\,\max\big\{\ax,\a\big\}-\a+\b\!\Big\}\,u\_1(\x).
\label{f}
\end{eqnarray}

Note that for every separably increasing $\fPx,$ function
$\fx\!:\Rk\to\R$ is well defined, i.e., the two terms in the right-hand
side of (\ref{f}) are finite. This follows from item (b) of
Proposition~\ref{P0}.

The main result of this section is

\begin{thm}
\label{T}
Suppose that $\fPx$ is a real-valued function defined on some
$P\subset\Rk$. Then $\fPx$ is monotonically extendible to $\Rk$ if and
only if $\fPx$ is separably increasing.
Moreover$,$ for every separably increasing $\fPx$, every strictly
increasing $u\_1(\x):\:\Rk\to\R$ that satisfies $(\ref{0<g1<1})$
and every $\a,\,\b\in\R\,\, s.t.\,\, \a<\b,$
the function $\fx$ defined by $(\ref{f})$ is a strictly increasing
extension of $\fPx$ to\/~$\Rk$.
\end{thm}

The same family of extensions $\fx$ can be generated using all strictly
increasing functions $u(\x)$ satisfying (\ref{ge}).

\begin{prop}
\label{P1}
For every separably increasing $\fPx$
and every strictly increasing $u(\x):$ $\Rk\to\R$ that
satisfies $(\ref{ge}),$ the function
\begin{eqnarray}
\fx =(\b-\a)^{-1} &\Big(&
\max\Big\{\ax-\a,\,\min\big\{\bx-\b,\,0\big\}\!\Big\}\big(\b-\ux\big)
\nonumber\\
&+&\min\Big\{\bx-\b,\,\max\big\{\ax-\a,\,0\big\}\!\Big\}\big(\ux-\a\big)
\Big)
\nonumber\\
&+&\ux
\label{f'}
\end{eqnarray}
is a strictly increasing extension of $\fPx$ to $\Rk$. This function
coincides with  $(\ref{f})$ where $u\_1(\x)$ is defined by $(\ref{u1})$.
\end{prop}

It is straightforward to verify that function (\ref{f'}) coincides with
(\ref{f}) where $u\_1(\x)$ is defined by~(\ref{u1}). This fact will be
employed in the proofs of Theorem~\ref{T} and some propositions, after
which the remaining statement of Proposition~\ref{P1} will follow from
Theorem~\ref{T}.

In (\ref{f'}) we assume $-\infty-\a=-\infty$ and $+\infty-\b=+\infty$
whenever $\ax=-\infty$ and $\bx=+\infty$, respectively, since $\a$ and
$\b$ are finite by definition.

To simplify (\ref{f'}), we partition $\Rk\setminus P$ into four
regions:
\begin{eqnarray}
\nonumber
A&=&\big\{\x\in \Rk\setminus P\;\big|\;
    (\exists\y\in P:\y<\x)\:\&\:(\exists\z\in P:\x<\z)\big\},\\
\nonumber
L&=&\big\{\x\in \Rk\setminus P\;\big|\;
    (\neg\exists\y\in P:\y<\x)\:\&\:(\exists\z\in P:\x<\z)\big\},\\
\nonumber
U&=&\big\{\x\in \Rk\setminus P\;\big|\;
    (\exists\y\in P:\y<\x)\:\&\:(\neg\exists\z\in P:\x<\z)\big\},\\
N&=&\big\{\x\in \Rk\setminus P\;\big|\;
     (\neg\exists\y\in P:\y<\x)\:\&\:(\neg\exists\z\in P:\x<\z)\big\}.
\label{ALUN}
\end{eqnarray}

Obviously,  every two of these regions have empty meet and
$\Rk=P\cup A\cup L\cup U\cup N$.

\begin{prop}
\label{P2}
For every separably increasing $\fPx,$ the
function $\fx$ defined by $(\ref{f'})$ can be represented as follows$:$
\begin{equation}
\fx=\cases{\fPx,                           &$\x\in P,$\cr
           \min\big\{\bx-\b,\,0\big\}+\ux, &$\x\in L,$\cr
           \max\big\{\ax-\a,\,0\big\}+\ux, &$\x\in U,$\cr
             \ux,                          &$\x\in N,$\cr
           \mbox{not simplified expression $(\ref{f'}),$}
                                           &$\x\in A$.\cr
}
\label{f''}
\end{equation}
\end{prop}

Proposition~\ref{P2} clarifies the role of $\ux$ in the definition
of~$\fx$. According to~(\ref{f''}), $\ux$ determines the rate
of growth of $\fx$ on $L$ and $U$, and $\fx=\ux$ on the set $N$ which
consists of $>$-neutral points with respect to the elements of~$P$.

Now let us give one more representation for~$\fx$, which can be used
in~(\ref{f''}) when $\x\in A$. Define four regions of other nature
in~$\Rk$:
\begin{eqnarray*}
S_1&=&\big\{\x\in \Rk\/\big|\;\bx-\ax\le\b-\a\big\},\cr
S_2&=&\big\{\x\in \Rk\/\big|\;\bx-\ax\ge\b-\a\mbox{ and }\,\bx\le\b\big\},\cr
S_3&=&\big\{\x\in \Rk\/\big|\;\bx-\ax\ge\b-\a\mbox{ and }\,\ax\ge\a\big\},\cr
S_4&=&\big\{\x\in \Rk\/\big|\;\ax\le\a\mbox{ and }\,\bx\ge\b\big\}.\cr
\end{eqnarray*}
It is easily seen that $\Rk=S_1\cup S_2\cup S_3\cup S_4$.

\begin{prop}
\label{P3}
For every separably increasing $\fPx,$ the function $\fx$ defined by
$(\ref{f'})$ or $(\ref{f}),$ with $u\_1(\x)$ and $\ux$ related by
$(\ref{u1}),$ can be represented as follows$:$
\begin{equation}
\fx=\cases{\ax\big(1-u\_1(\x)\big)+\bx\,u\_1(\x), &$\x\in S_1,$\cr
             \bx+\ux-\b,                          &$\x\in S_2,$\cr
             \ax+\ux-\a,                          &$\x\in S_3,$\cr
             \ux,                                 &$\x\in S_4.$\cr
}
\label{f'''}
\end{equation}

\end{prop}

The regions $S_1,S_2,S_3,$ and $S_4$ meet on some parts of the border
sets $\bx-\ax=\b-\a$, $\ax=\a$, and $\bx=\b$. Accordingly, the
expressions of $\fx$ given by Proposition~\ref{P3} are concordant on
these intersections.

\section{Extendibility of arbitrary functions defined on Pareto sets}
\label{SPar}

Consider the case where $P$ is a Pareto set.

\begin{defn}
A set $P\subset\Rk$ is a {\em Pareto set\/} in $\Rk$ if there are no
$\x,\x'\in P$ such that $\x'>\x$.
\end{defn}

Observe that for every function $\fPx$ defined on a Pareto set $P,$ the
set $A$ introduced in~(\ref{ALUN}) is empty. It turns out that
such a function $\fPx$ is separably increasing if and only if it is
upper-bounded on lower sets and lower-bounded on upper sets (see
Definition~\ref{def23}). Based on this, the following corollary from
Theorem~\ref{T} and Proposition~\ref{P2} is true.

\begin{cor}
\label{cor}
Suppose that $\fPx:\: P\to\R$ is a function defined on a Pareto set
$P\subset\Rk$. Then $\fPx$ can be strictly monotonically extended
to $\Rk$ if and only if $\fPx$ is
upper-bounded on lower sets and
lower-bounded on upper sets.
Moreover$,$ for such a function $\fPx$ and every strictly increasing
function $\ux\!:\:\Rk\to\R$ that satisfies $(\ref{ge}),$ the function
\begin{equation}
\fx=\cases{\fPx,                              &$\x\in P,$\cr
           \min\big\{\bx,\,\b\big\}-(\b-\ux), &$\x\in L,$\cr
           \max\big\{\ax,\a\big\}+\ux-\a,     &$\x\in U,$\cr
           \ux,                               &$\x\in N,$ \cr
}
\label{Pare_f}
\end{equation}
provides a monotone extension of $\fP(\x)$ to~$\Rk$.
\end{cor}

On regions $L$ and $U$, $f(\x)$ is expressed through the ``relative''
functions $\b-\ux$ and $\ux-\a$; this is equivalent to the form given
in~(\ref{f''}). A result closely related to Corollary~\ref{cor} was
used in Chebotarev and Shamis~(1998) to construct an implicit form of
monotonic scoring procedures for preference aggregation.

\section{The extension problem in the context of utility theory}
\label{SUti}

Extensions and utility representations of partial orders were studied
since Zorn's lemma and the Szpilrajn theorem according to which every
strict partial order extends to a strict linear order.

In general, neither  strict partial orders nor their linear extensions
must have utility representations. Let us discuss connections between
these extensions and utility representations in more detail. Recall that

\begin{defn}
A {\em utility representation\/} of a strict
partial order $\succ$ on a set $X$ is a function $u:\,\,X\to\R$ such
that for every\footnote{For uniformity, we designate the elements of
$X$ by boldface letters.}
$\x,\y\in X,$
\begin{eqnarray}
\x\succ \y&\Rightarrow& u(\x)>u(\y),
\label{succ}
\\
\x\approx \y&\Rightarrow& u(\x)=u(\y),
\label{approx}
\end{eqnarray}
where$,$ by definition$,$
\begin{eqnarray*}
\x\approx \y &\Leftrightarrow&
     [\,\forall \z\in X,\; \x\sim \z\Leftrightarrow \y\sim \z\,],\\
\x\sim \y &\Leftrightarrow&
     [\,\x\not\succ \y\mbox{ and }\y\not\succ \x\,].
\end{eqnarray*}
\end{defn}

A sufficient condition for a strict partial order $\succ$ to have
a utility representation is the existence of a countable and dense
(w.r.t.\ the induced strict partial order) subset in the factor
set $X/\!\approx$ (see Debreu, 1964; Fishburn, 1970). Generally, this
sufficient condition is not necessary, but if $\succ$ is a strict
weak order, then it is necessary.

The Paretian $>$ relation on $\Rk$ is a special strict partial order.
Any lexicographic order on $\Rk$ is its strict linear extension. Such
extensions have no utility representations, whereas the $>$ relation
has a wide class of utility representations.\footnote{They do exist,
say, because $\Rk$ contains countable and dense (w.r.t.\ $>$) subsets,
one of which being the set of vectors with rational coordinates.}
These are all functions strictly increasing in all coordinates.

Every such a strictly increasing
function induces a strict weak order on $\Rk$ that
extends~$>$. Naturally, not all strict weak orders that extend $>$
can be obtained in this manner. A sufficient condition for
such representability is
the Archimedean property which ensures the existence of a
countable and dense (w.r.t.\ the strict weak order) subset
in~$\Rk$.

Thus, the utility representations of $>$ induce a special class of
strict weak orders that extend~$>$. Such a strict weak order determines
its utility representation up to arbitrary monotone transformations
(some related specific results are given in Morkeli\=unas, 1986b).

In the previous sections, we considered
a (utility) function $\fP$ that is defined on a subset $P$ of $\Rk$ and
represents the restriction of $>$ to~$P$. If $P$ is a Pareto subset,
then this imposes no constraints on~$\fP$. The problem was to find
conditions under which there exist functions $f$ that $(\ast)$~reduce to
$\fP$ on $P$ and $(\ast\ast)$~represent $>$ on $\Rk$, and to provide
an explicit form of such functions.

Observe that every strictly increasing function $\fP$ induces some
strict weak order $\succ$ on $P$ that contains the restriction of $>$
to~$P$. For $\succ$, there exists a countable and dense
(w.r.t.\ $\succ$) subset in the factor set $P/\!\sim$, where $\sim$
corresponds to~$\succ$. Combining this subset with the set of vectors
in $\Rk$ that have rational coordinates, we obtain a countable and dense
subset in the factor set $\Rk/\!\sim$ (the $\approx$ relation that
corresponds to $>$ is the identify relation). Consequently, there exist
utility functions $g$ that represent $T(>\cup\succ)$ on $\Rk$, where
$T(\cdot)$ designates transitive closure.

The restriction of such a function $g$ to $P$ need not coincide
with $\fP$, but it is related with $\fP$ by a strictly increasing
transformation, since they represent the same weak order~$\succ$. Thus,
we obtain

\begin{prop}
\label{Put}
For every strictly increasing function $\fP$ on $P\subset\Rk,$ there
exists a strictly increasing map $\varphi$ of the range of $\fP$ to $\R$
such that $\varphi(\fP)$ is monotonically extendible to~$\Rk$. If $\fP$
is not strictly increasing$,$ then there are no such maps.
\end{prop}

Proposition~\ref{Put} elucidates a difference between separably
increasing and strictly increasing functions $\fP$ with respect to the
extendibility. The former functions are monotonically extendible
to $\Rk$ (Theorem~\ref{T}), whereas for the latter functions, only some
their strictly increasing transformations are extendible in the general
case.

\section{Extending utility functions on arbitrary sets}
\label{abstSet}

A shortcoming of the definition~(\ref{f}) is that it provides
discontinuous extensions even for continuous functions~$\fP$. On
the other hand, the extension technique does not strongly rely on the
structure of~$\Rk$. Indeed, the connection of $\fx$ with $\Rk$, its
domain, reduces to the dependence on $\ax,\bx,$ and $\ux$, which have a
rather general nature. This enables one to translate the above results
to an abstract set $X$ substituted for $\Rk$ and any strict partial
order $>$ on $X$ substituted for the Paretian $>$ relation on~$\Rk$.
The author thanks Andrey Vladimirov for his suggestion to consider the
problem in this general setting. To formulate a counterpart of
Theorem~\ref{T}, we will use the following notation:

$X$ is a nonempty set;\\
$>$ is a fixed strict partial order on $X$;\\
$P\subset X;\,$ $>\!\_P$ is the restriction of $>$ to~$P$;\\
$\z\in X$ is a {\em maximal element\/} of $>$ iff $\x>\z$ for no
$\x\in X$;\\
$\y\in X$ is a {\em minimal element\/} of $>$  iff $\x<\y$ for no $\x\in X$.

We will also explore all the above definitions and formulas (except
for~(\ref{arctan}) and Definition~\ref{def1}), where the
substitution of $X$ for $\Rk$ is implied.

For arbitrary strict partial orders---which may have maximal and minimal
elements---a somewhat stronger condition should be used in the
definition of separably increasing functions. Let
$$
\widetilde X=X\cup\,\{-\infty,+\infty\}.
$$

Extend $>$ to $\widetilde X$ in the usual way to get a strict partial
order on $\widetilde X$: $+\infty>-\infty$ and $+\infty>\x>-\infty$ for
every $\x\in X$. We use the same symbol $>$ for a strict partial order
on $X$, for its extension to $\widetilde X$, and for the ordinary
``greater'' relation on $\R$, which should not lead to confusion.

\begin{defn}
\label{def1'}
A~function $\fP(\x)$ defined on $P\subset X$ is {\em separably
increasing\/} if for any $\x,\x'\in\widetilde X,\:$ $\x'>\x$ implies
$b(\x')>a(\x)$.
\end{defn}

Since for every $\fP$, $b(+\infty)=+\infty$ and $a(-\infty)=-\infty,$
Definition~\ref{def1'} implies that for a separably increasing function,
$a(\x)<+\infty$ and $b(\x)>-\infty$ whenever $x\in X$ is a maximal and a
minimal element of $>$, respectively. If $>$ has neither maximal nor
minimal elements (as for the Paretian $>$ relation on~$\Rk$), then the
replacement of $\widetilde X$ with $X$ in Definition~\ref{def1'} does
not alter the class of separably increasing functions.

Obviously, if $>$ has a utility representation, then it has bounded
utility representations (which can be constructed, say, by
transformations like~(\ref{arctan})).

\begin{thm}
\label{T'}
Suppose that a strict partial order $>$ defined on $X$ has a utility
representation$,$ $P\subset X,$ and $\fP\!: P\to\R$ is a utility
representation of~$>\!\_P$. Then $\fP$ is monotonically extendible to
$X$ if and only if $\fP$ is separably increasing.\\
Moreover$,$ if $u\_1(\x)$ is a utility representation of $>$
that satisfies~$(\ref{0<g1<1})$ and $\a<\b,$
then $(\ref{f})$ provides a monotone extension of $\fP$
to~$X$.
\end{thm}

If a function $\fP$ is not separably increasing according to
Definition~\ref{def1'} but it satisfies this definition with
$\widetilde X$ replaced by $X$, then (\ref{f}) can be used to obtain
extensions of $\fP$ in terms of ``quasiutility'' functions $f\!:
X\to\widetilde{\R}$.

Propositions~\ref{P0},~\ref{P1},~\ref{P2}, and~\ref{P3} are preserved
in the case of arbitrary~$X$.
Proposition~\ref{Put}, as well as Theorem~\ref{T'}, is valid for the
strict partial orders $>$ on $X$ that have a utility representation.
The proofs given in Section~\ref{proofs} are conducted for this general
case.

Interesting further problems are describing the complete class of
extensions of $\fP$ to~$\Rk$ (and to $X$), specifying necessary and
sufficient conditions for the existence of {\em continuous\/}
extensions and constructing them, and considering the extension problem
with an ordered extension of the field $\R$ as the range of $\fP,u,$
and~$f$. An interesting result on the existence of continuous utility
representations for weak orders on $\Rk$ was obtained in
Morkeli\=unas~(1986a).

\section{Proofs}
\label{proofs}

The proofs of all propositions and Corollary~\ref{cor} are given in the
general case of a strict partial order on an abstract set $X$ (see
Section~\ref{abstSet}, especially, Definition~\ref{def1'}).
Theorem~\ref{T'} which generalizes Theorem~\ref{T} is proved
separately.

\begin{pf*}{Proof of Proposition~\ref{P0}.}
(a) Assume that $\fPx$ is not strictly increasing. Then there are
$\x,\x'\in P$ such that $\x'>\x$ and $\fP(\x')\le\fPx$. Then, by
(\ref{a>f>b}),  $b(\x')\le\fP(\x')\le\fPx\le\ax$ holds$,$ i.e., $\fPx$
is not separably increasing.

(b) Let $P'$ be a lower set. Then there exists $\aa\in X$ such that
$P'=\big\{\x\mid\x\le\aa,\x\in P\big\}$. Consider any
$\aa'\in\widetilde X$ such that $\aa'>\aa$. Since $\fPx$ is separably
increasing, $b(\aa')>a(\aa)$. Therefore, $a(\aa)<+\infty$. Since
$a(\aa)=\sup\big\{\fP(\y)\mid\y\in P'\big\}$, $\fPx$ is upper-bounded
on~$P'$. Similarly, $\fPx$ is lower-bounded on upper sets.

(c) Let $\x\in X.$ By (a), if $\y,\z\in P,\,$ $\z\ge\x,$ and
$\y\le\x,$ then $\fP(\z)\ge\fP(\y)$. Having in mind that
$\sup\emptyset=-\infty$ and $\inf\emptyset=+\infty$, we
obtain~$b(\x)\ge a(\x)$.

(d) By (a), for every $\x\in P,\,$ $b(\x)=\fPx$ and $a(\x)=\fPx$ hold.
This completes the proof.
\end{pf*}

Now we prove Proposition~\ref{P3}; thereafter Proposition~\ref{P3} will
be used to prove Proposition~\ref{P2} and Theorem~\ref{T}.

\begin{pf*}{Proof of Proposition~\ref{P3}.}
Let $\x\in S_1.$ Since $\bx-\ax\le\b-\a,$ we have
$$
\min\big\{\bx,\,\b\big\}-\ax\le\b-\a,
$$
$$
\bx-\max\big\{\ax,\a\big\}\le\b-\a,
$$
whence
$$
\ax\ge\min\big\{\bx,\,\b\big\}-\b+\a,
$$
$$
\bx\le\max\big\{\ax,\a\big\}-\a+\b.
$$
Therefore, (\ref{f'}) reduces to $\fx=\ax\big(1-u\_1(\x)\big)+\bx u\_1(\x)$.

Let $\x\in S_2.$ Inequalities $\bx-\ax\ge\b-\a$ and $\bx\le\b$ imply
$\ax\le\a$, therefore, (\ref{f}) reduces to $\fx=\bx+\ux-\b$.

Let $\x\in S_3.$ Inequalities $\bx-\ax\ge\b-\a$ and $\ax\ge\a$ imply
$\bx\ge\b$, therefore, (\ref{f}) reduces to $\fx=\ax+\ux-\a$.

The proof for the case $\x\in S_4$ is straightforward.
\end{pf*}

\begin{pf*}{Proof of Proposition~\ref{P2}.}
Let $\x\in P.$ Then, by item (d) of Proposition~\ref{P0},
$\bx=\ax=\fPx$, therefore, by (\ref{ge}), $\bx-\ax\le\b-a$, hence
$\x\in S_1$. Using Proposition~\ref{P3}, we have
$\fx=\fPx\big(1-u\_1(\x)\big)+\fPx\, u\_1(\x)=\fPx$.

Let $\x\in U$. Then
$\bx=+\infty$, hence (\ref{f'}) reduces to
$\fx=\max\big\{\ax-\a,\,0\big\}+\ux$. Similarly,
if  $\x\in L$, then
$\ax=-\infty$ and (\ref{f'}) reduces to
$\fx=\min\big\{\bx-\b,\,0\big\}+\ux$.

Finally, if $\x\in N,$ then $\ax=-\infty$ and $\bx=+\infty$, whence
$\ax<\a$ and $\bx>\b$, and Proposition~\ref{P3} provides $\fx=\ux$.
\end{pf*}

\begin{pf*}{Proof of Theorem~\ref{T}.}
Let $\fPx$ be separably increasing.
By Proposition~\ref{P2}, the restriction of $\fx$ to $P$
coincides with~$\fPx$.

Prove that $\fx$ is strictly increasing on~$\Rk$. This can be
demonstrated directly by analyzing equation~(\ref{f}). Here, we
give another proof, which does not require any additional
calculations with $\,\min\,$ and $\,\max\,$.

By Proposition~\ref{P3}, function~(\ref{f}) coincides with
(\ref{f'''}), where $\ux$ is related with $u\_1(\x)$ by~(\ref{u1}).

Suppose that $\x,\x'\in\Rk$ and $\x'>\x$. Then, by (\ref{abmono}) and
the strict monotonicity of $\ux$ and $u\_1(\x)$, we have
\begin{eqnarray}
u(\x')   &>&\ux,     \cr
u\_1(\x')&>&u\_1(\x),\cr
a(\x')  &\ge&\ax,    \cr
b(\x')  &\ge&\bx.
\label{pre}
\end{eqnarray}

Suppose first that $\x$ and $\x'$ belong to the same region:
$S_2,S_3,$ or $S_4$. Then (\ref{pre}) yields
\begin{eqnarray}
b(\x')+u(\x')-\b&>&\bx+\ux-\b, \cr
a(\x')+u(\x')-\a&>&\ax+\ux-\a, \cr
          u(\x')&>&\ux,
\label{mo}
\end{eqnarray}
hence, by (\ref{f'''}), $f(\x)$ is strictly increasing on each of these
regions.

If $\x,\x'\in S_1$, then by (\ref{f'''}), (\ref{pre}), (\ref{0<g1<1}),
and (c)~of Proposition~\ref{P0},
\begin{eqnarray}
f(\x')-\fx
&\ge& a(\x)\big(1-u\_1(\x')\big)+b(\x)u\_1(\x')\nonumber\\
&-&   a(\x)\big(1-u\_1(\x )\big)-b(\x)u\_1(\x )\nonumber\\
&=&   \big(b(\x)-a(\x)\big)\big(u\_1(\x')-u\_1(\x)\big)\ge0.
\end{eqnarray}

This implies that $f(\x')=\fx$ is possible only if $b(\x')=b(\x)$
and $b(\x)=a(\x)$, i.e., only if $b(\x')=a(\x)$. The last equality
is impossible, since $\fPx$ is separably increasing by assumption.
Therefore, $f(\x')>\fx$, and $\fx$ is strictly increasing
on~$S_1$.

Let now $\x$ and $\x'$ belong to different regions $S_i$ and~$S_j$.
Consider the points that correspond to $\x$ and $\x'$ in the
3-dimensional space with coordinate axes $a(\cdot)$, $b(\cdot)$, and
$u(\cdot)$ and connect these two points, $\big(\ax,\bx,\ux\big)$ and
$\big(a(\x'),b(\x'),u(\x')\big)$, by a line segment.
The projection of the line segment and the borders of the regions
$S_1,S_2,S_3,$ and $S_4$ to the plane $u=0$ are illustrated
in Fig.~\ref{F1}.

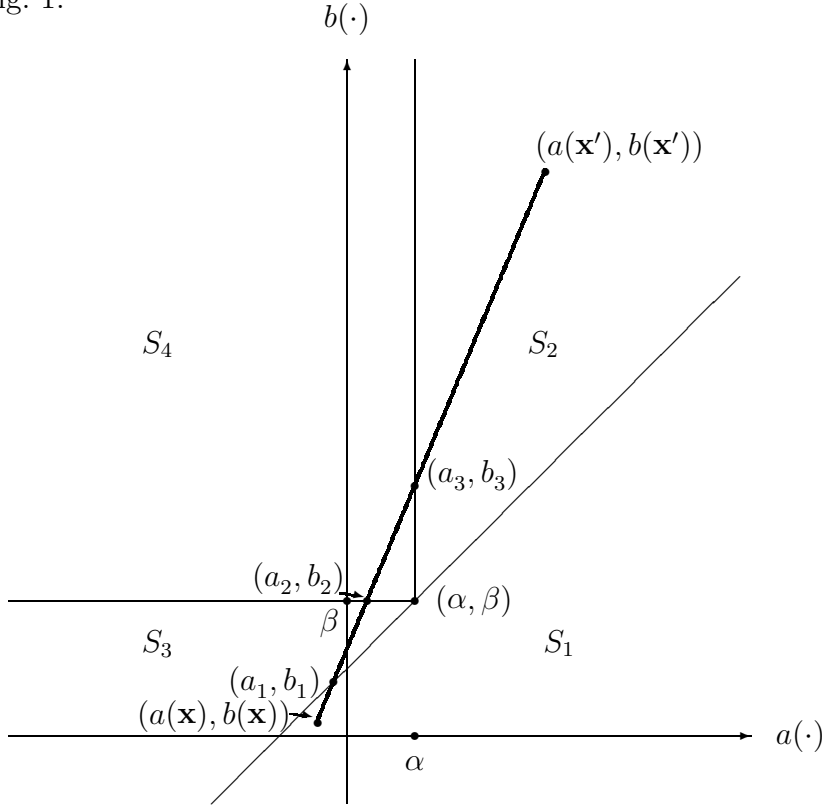
\begin{figure}[ht]
\unitlength 0.90mm
\linethickness{0.3pt}
\begin{picture}(133.33,123.67)
\put(60.00,10.00){\vector(0,1){110}}
\put(10.00,20.00){\vector(1,0){109.67}}
\linethickness{0.6pt}
\put(40.00,10.00){\line(1,1){78.00}}
\put(10.00,40.00){\line(1,0){60.00}}
\put(70.00,40.00){\line(0,1){80.00}}
\linethickness{1.2pt}
\multiput(55.67,22.33)(0.12,0.29){281}{\line(0,1){0.29}}
\linethickness{0.6pt}
\put(89.33,103.45){\circle*{1.20}}
\put(55.67,22.00){\circle*{1.20}}
\put(70.00,20.00){\circle*{1.20}}
\put(60.00,40.00){\circle*{1.20}}
\put(58.00,28.00){\circle*{1.20}}
\put(63.00,40.00){\circle*{1.20}}
\put(70.00,57.00){\circle*{1.20}}
\put(70.00,40.00){\circle*{1.20}}
\put(123.33,20.00){\makebox(0,0)[lc]{$a(\cdot)$}}
\put(60.00,123.67){\makebox(0,0)[cb]{$b(\cdot)$}}
\put(73.00,40.00){\makebox(0,0)[lc]{$(\a,\b)$}}
\put(51.67,23.00){\makebox(0,0)[rc]{$\big(\ax,\bx\big)$}}
\put(56.14,28.00){\makebox(0,0)[rc]{$(a\_1,b\_1)$}}
\put(59.67,41.00){\makebox(0,0)[rb]{$(a\_2,b\_2)$}}
\put(85.33,56.33){\makebox(0,0)[rb]{$(a\_3,b\_3)$}}
\put(112.67,104.67){\makebox(0,0)[rb]{$\big(a(\x'),b(\x')\big)$}}
\put(70.00,17.00){\makebox(0,0)[ct]{$\a$}}
\put(58.80,38.80){\makebox(0,0)[rt]{$\b$}}
\put(89.00,35.67){\makebox(0,0)[lt]{$S_1$}}
\put(89.00,76.00){\makebox(0,0)[cb]{$S_2$}}
\put(32.00,31.67){\makebox(0,0)[cb]{$S_3$}}
\put(32.00,76.00){\makebox(0,0)[cb]{$S_4$}}
\put(54.73,22.85){\vector(4,-1){0.2}}
\multiput(51.57,23.30)(0.27,-0.055){8}{\line(1,0){0.54}}
\put(62.29,40.66){\vector(4,-1){0.2}}
\multiput(58.95,41.11)(0.29,-0.055){8}{\line(1,0){0.59}}
\end{picture}
\caption{An example of line segment
$\big[\big(\ax,\bx,\ux\big),\big(a(\x'),b(\x'),u(\x')\big)\big]$ in the
space with coordinate axes $a(\cdot)$, $b(\cdot)$, and $u(\cdot)$
projected to the plane~$u=0$.
}
\label{F1}
\end{figure}

Suppose that $(a\_1,b\_1,u\_1),\ldots,(a\_p,b\_p,u\_p)$, $p\le3$, are
the points where the line segment between $\big(\ax,\bx,\ux\big)$ and
$\big(a(\x'),b(\x'),u(\x')\big)$ crosses the borders of the regions.
Then

\begin{equation}
\ax\le a\_1\le\cdots\le a\_p\le a(\x'),
\label{a_le}
\end{equation}
\begin{equation}
\bx\le b\_1\le\cdots\le b\_p\le b(\x'),
\label{b_le}
\end{equation}
\begin{equation}
\ux<u\_1<\cdots<u\_p<u(\x')
\label{u_le}
\end{equation}
with strict inequalities in (\ref{a_le}) or in (\ref{b_le}) (or in the
both).

Consider $\fx$ represented by (\ref{f'''}) as a function $\fhat(a,b,u)$
of $\ax,\bx$, and $\ux$. Then, using the fact that $\fhat(a,b,u)$ is
nondecreasing in $a$ and $b$ on each region, strictly increasing in $u$
on $S_2, S_3,$ and $S_4$, and strictly increasing in $u$ on $S_1$ unless
$\ax=b(\x')$ (which is not the case, since $\fPx$ is separably increasing),
we obtain
\begin{eqnarray}
\fx
&=&\fhat\big(\ax,\bx,\ux\big)<\fhat(a\_1,b\_1,u\_1)<\cdots
                             <\fhat(a\_p,b\_p,u\_p)\cr
&<&\fhat\big(a(\x'),b(\x'),u(\x')\big)=f(\x').
\label{fx<fx'}
\end{eqnarray}

This completes the proof that $\fx$ is strictly increasing.

It remains to prove that if $\fPx$ is not separably increasing,
then it cannot be strictly monotonically extended to~$\Rk$.
Indeed, if there are $\x,\x'\in\Rk$ such that $\x'>\x$
and $b(\x')\le\ax$, then strict monotonicity of $\fx$ requires
$f(\x')\le b(\x')$ and $\fx\ge\ax$ to hold, whence $\fx\ge f(\x')$,
and strict monotonicity is violated. Theorem~\ref{T} is proved.
\end{pf*}

\begin{pf*}{Proof of Corollary~\ref{cor}.}
Suppose that a strict partial order $>$ on $X$ has a utility
representation.
Let $\fP$ be a utility
representation of $>\!\_P$, where $>\!\_P$ is the restriction of $>$ to a
Pareto set $P\subset X.$ Then, since the $>$ relation is
transitive, the set $A$ is empty. It remains to prove that $\fP$ is
separably increasing if and only if it is upper-bounded on lower sets
and lower-bounded on upper sets. If $\fP$ is separably increasing,
then these boundedness conditions are satisfied by
Proposition~\ref{P0}. Suppose now that $\fP$ is upper-bounded on
lower sets and lower-bounded on upper sets and assume that $\fP$ is
not separably increasing. Then there exist $\x,\x'\in\widetilde X$ such
that $\x'>\x$ and $b(\x')\le\ax$. This is possible only if
(a)~$b(\x')=+\infty$ or (b)~$\ax=-\infty$ or (c)~there are $\y,\z\in P$
such that $\y\le\x$ and $\z\ge\x'.$ However, in (a), $\ax=+\infty$ and
$\x\in X$, hence  $\fP$ is not upper-bounded on a lower set; in (b),
$b(\x')=-\infty$ and $\x'\in X$, hence  $\fP$ is not lower-bounded on
an upper set; in (c), by transitivity, $\z>\y,$ which contradicts the
definition of Pareto set. Therefore, $\fP$ is separably increasing,
and the corollary is proved.
\end{pf*}

\begin{pf*}{Proof of Theorem~\ref{T'}.}
If $\fP$ is not separably increasing,
then it is not monotonically extendible to~$X$. Indeed, if the
implication $\x'>\x\Rightarrow b(\x')>a(\x)$ is violated for some
$\x,\x'\in X$, then the argument is the same as for Theorem~\ref{T}.
Otherwise, if this implication is violated for $\x\in\widetilde
X\setminus X$, then $\x'>\x$ implies $\x=-\infty$, and since $b(\x')\le
a(\x)=-\infty$, $\x'\in X$ holds, hence $f(\x')$ cannot be evaluated
without violation of~(\ref{succ}); the case of $\x'\in\widetilde
X\setminus X$ is considered similarly.

The argument in the proof of Theorem~\ref{T} works here with no change
to demonstrate that if $>$ has a utility representation and $\fP$ is
separably increasing, then every $f$ defined by (\ref{f}) satisfies
condition~(\ref{succ}) of utility representability. It remains to show
that $f$ satisfies~(\ref{approx}). Let $\x,\y\in X$ and $\x\approx \y$.
Then for every $\z\in X$, $[\x\succ \z\Leftrightarrow \y\succ \z$ and
$\z\succ \x\Leftrightarrow \z\succ \y]$ (see, e.g., Fishburn, 1970).
Consequently, $a(\x)=a(\y)$ and $b(\x)=b(\y)$. Since $u$ represents
$>$, $u(\x)=u(\y)$ holds. Therefore, $f(\x)=f(\y)$. This completes the
proof.
\end{pf*}
\bigskip

The author thanks Elena Yanovskaya, Elena Shamis, Andrey
Vladimi\-rov, and Fuad Aleskerov for helpful discussions.

\section*{References}

\begin{description}

\item[]
{Chebotarev, P.Yu. and E. Shamis,} 1998,
Characterizations of scoring methods for preference aggregation,
{Annals of Operations Research}
80, 299--332.

\item[]
{Debreu, H.,} 1964,
Continuity properties of Paretian utility,
{International Economic Review}
5, 285--293.

\item[]
{Fishburn, P.,} 1970,
{Utility Theory and Decision Making}
(Wiley, New York).

\item[]
{Morkeli\=unas, A.,} 1986a, On the existence of a
continuous superutility function,
{Lietuvos Matematikos Rinkinys / Litovskiy Matematicheskiy Sbornik}
26, 292--297
(Russian).

\item[]
{Morkeli\=unas, A.,} 1986b,
On strictly increasing numerical transformations and the Pareto
condition,
{Lietuvos Matematikos Rinkinys / Litovskiy Matematicheskiy Sbornik}
26, 729--737
(Russian).

\end{description}

\end{document}